\documentclass[11pt]{article}

\input{amssym.def}
\input{amssym}
\usepackage{eufrak}

\newtheorem{theorem}{Theorem}

\newtheorem{conjecture}[theorem]{Conjecture}
\newtheorem{corollary}[theorem]{Corollary}
\newtheorem{definition}[theorem]{Definition}

\newtheorem{proposition}[theorem]{Proposition}
\newtheorem{remark}[theorem]{Remark}

\def\h{h}
\def\hh{\hbar}
\def\pa{{\partial}}
\def\hl{\hat{l}}
\def\hL{\hat{L}}

\def\hatt{\hat{\Theta}}

\def\un{\underline}
\def\qq{q^{-1}}

\def\Tr{\mathrm{Tr}}
\def\Trr{\Tr_R}
\def\rh{r_\hh}

\def\hLL{\hat{{\cal L}}}
\def\LL{{\cal{L}}}
\def\MM{{\cal{M}}}

\def\DD{{\cal D}}
\def\NN{{\cal N}}

\def\si{{\sigma}}
\def\al{{\alpha}}

\def\hpa{\hat{\pa}}
\def\hD{\hat{D}}
\def\de{\delta}
\def\De{\Delta}
\def\ot{\otimes}
\def\C{{\Bbb C}}

\def\vv{V^{\otimes 2}}

\def\ov{\overline}

\def\la{{\lambda}}

\def\be{\begin{equation}}
\def\ee{\end{equation}}

\begin{document}

\makeatletter
\renewcommand{\theequation}{{\thesection}.{\arabic{equation}}}
\@addtoreset{equation}{section}
\makeatother

\title{\rightline{\normalsize\it To the memory of A.M. Vinogradov}
\bigskip
Quantum vector fields via quantum doubles and their applications}
\author{\rule{0pt}{7mm} Dimitry Gurevich\thanks{gurevich@ihes.fr}\\
{\small\it Institute for Information Transmission Problems}\\
{\small\it Bolshoy Karetny per.~19,  Moscow 127051, Russian Federation}\\
\rule{0pt}{7mm} Pavel Saponov\thanks{Pavel.Saponov@ihep.ru}\\
{\small\it
National Research University Higher School of Economics,}\\
{\small\it 20 Myasnitskaya Ulitsa, Moscow 101000, Russian Federation}\\
{\small \it and}\\
{\small \it
Institute for High Energy Physics, NRC "Kurchatov Institute"}\\
{\small \it Protvino 142281, Russian Federation }}
\maketitle

\begin{abstract}
By treating generators of the reflection equation algebra corresponding to a Hecke symmetry as quantum analogs of  vector fields, we exhibit the corresponding Leibniz rule
via the so-called quantum doubles. The role of the function algebra in such a double is attributed to another copy of the reflection equation algebra. We consider two types
of quantum doubles: these giving rise to the quantum analogs of left vector fields acting on the function algebra and those giving rise to quantum analogs of the adjoint vector
fields acting on the same algebra. Also, we introduce quantum partial derivatives in the generators of the reflection equation algebra and then at the limit $q\rightarrow 1$ we get
quantum partial derivatives on the enveloping algebra $U(gl_N)$ as well as on a certain its extension.
\end{abstract}

\maketitle

\section{Introduction}

The central problem of any attempt to generalize the notion of a vector field on a manifold is the definition of the corresponding generalized Leibniz rule. Thus, the Leibniz rule for 
su\-per-fields acting on a su\-per-va\-ri\-ety should take into account  the parity of all elements involved. Such su\-per-fields generate super-Lie algebras. In the late 80's one of the 
authors introduced analogs of Lie algebras, associated with involutive symmetries, that is braidings\footnote{Recall that by a braiding we mean a linear  operator $R:\vv\to \vv$ 
subject to the so-called braid relation
$$
(R\ot I)(I\ot R)(R\ot I)=(I\ot R)(R\ot I)(I\ot R).
$$ } 
$R:\vv\to\vv$, such that $R^2=I\otimes I$ (see \cite{G} and references therein). Hereafter, $V$ is a vector space of dimension $N$ over the field $\C$ and $I$ is the identity 
operator in $V$ or the $N\times N$ unit matrix. The Leibniz rule for the corresponding vector fields was also introduced via such symmetries.

Latter it became clear that many aspects of this approach can be formulated in terms of the so-called reflection equation (RE) algebras, associated with Hecke symmetries,
that is  braidings $R$ such that
\be
R^2 = I\otimes I +(q-q^{-1})R,\qquad q\in \C\setminus \{\pm 1,0\}.
\label{H-cond}
\ee

With a given Hecke symmetry $R$ we associate two forms of the RE algebra: the non-modified $\LL(R)$ and modified $\hLL(R)$ ones. They are respectively generated by
entries of the matrices $L=(l_i^j)_{1\leq i,j \leq N}$ and $\hL=(\hl_i^j)_{1\leq i,j \leq N}$, which are subject to the systems of relations
\be
R_{12}L_1R_{12} L_1-L_1R_{12} L_1R_{12}=0,
\label{RE}
\ee
\be
R_{12}\hL_1R_{12} \hL_1-\hL_1R_{12}\hL_1R_{12}=R_{12}\hL_1-\hL_1 R_{12},
\label{mRE}
\ee
where $L_1 = L\otimes I$ and similarly for all other matrices.

Note that these algebras are isomorphic to each other. Their isomorphism can be established by the following relation on the {\em generating matrices} $L$ and $\hL$:
\be
L=I-(q-\qq)\, \hL.
\label{iso}
\ee

Thus, in fact, we are dealing with one algebra realized in two different sets of generators: $\{l_i^j\}$ and $\{\hl_i^j\}$. Observe that the isomorphism (\ref{iso}) fails at $q=\pm 1$. If a Hecke
symmetry $R$ is a deformation of the usual flip  $P$, then at the limit $q\to 1$ the algebra $\LL(R)$ tends to $\mathrm{Sym}(gl_N)$, whereas $\hLL(R)$ tends to the algebra $U(gl_N)$.

The generators of the RE algebra (modified or not) are usually considered as quantum (or braided) analogs of the classical vector fields. Habitually, the role of the function algebra
is played by the corresponding  RTT algebra (see \cite{IP} and references therein). In \cite{GPS2, GPS3} we suggested a different version of the quantum calculus, in which the role of the
function algebra was played by another copy of the RE algebra or by an algebra, generated by the space $V$ or its dual $V^*$.

All these examples can be included into the family of the so-called quantum doubles. By a quantum double we mean a couple of associative algebras $(A,B)$ endowed with a permutation  map $\si:A\ot B\to B\ot A$, satisfying a set of certain axioms (see the next section). On the base of this map one can introduce the structure of a unital associative algebra on the linear space $B\otimes A$. If the algebra $A$ is endowed with an algebra homomorphism $\varepsilon_A: A\to \C$, then it is possible to define an action of the algebra $A$ onto $B$ (see 
formula (\ref{act}) below) and this action will be compatible with the algebraic structures of both components.

The objective of this article is two-folded. First, we  compare two quantum doubles, giving rise to quantum analogs of the left vector fields and to quantum analogs of the adjoint vector
fields. These quantum vector fields act on the RE algebra, which plays the role of the function algebra $\mathrm{Sym}(gl_N)$. This construction is exhibited in Section 2 in terms of  quantum doubles. 

Second, we discuss some applications of these two classes of quantum vector fields. Namely, in Section 3 we consider invariant differential operators, corresponding to the central elements of the RE algebra and represented by the quantum left vector fields, and formulate a conjecture on their spectrum.

In Section 4 we establish a series of matrix Capelli identities which are quantum generalizations of results by A.Okounkov  (see \cite{Ok1,Ok2}). These generalizations  are based on the quantum analogs of the partial derivatives in the generators of an RE algebra, defined in Section 2. Moreover, in Section 6 we get analogs of the partial derivatives on the enveloping algebra
$U(gl_N)$ by passing to the limit $q\rightarrow 1$ in the corresponding structures of the RE algebra. As for the quantum adjoint vector fields, we show that they can be reduced to 
the so-called ``quantum orbits", which are some quotients of the RE algebras. This is done in Section 5.

\section{Quantum doubles and partial derivatives}

\begin{definition}
A couple $(A, B)$ of unital associative algebras $A$ and $B$ is called {\it an associative double} if it is equipped with a linear invertible {\it permutation map}
$$
\sigma :A \ot B\to B\ot A,
$$
which meets the following requirements:
$$
\sigma\circ(\mu_{A}\otimes\mathrm{id}_B) =(\mathrm{id}_B\otimes\mu_{A})\circ\sigma_{12}\circ\sigma_{23} \quad\mathrm{on}\quad A\ot A\ot B,
$$
$$
\sigma\circ(\mathrm{id}_A\otimes\mu_{B}) =(\mu_{B}\otimes \mathrm{id}_A)\circ\sigma_{23}\circ\sigma_{12} \quad\mathrm{on}\quad A\ot B\ot B,
$$
$$
\sigma(1_A\ot b)=b\ot 1_A,\quad \sigma(a\ot 1_B)=1_B\ot a \qquad\forall\, a\in A,\, \forall\,b\in B.
$$
Here $\mu_A: A\ot A\to A $ is the multiplication in the algebra $A$, $1_A$ is its unit element, and similarly for $B$.
\end{definition}

In general, if the algebras $A$ and $B$ are introduced via  relations  on their  generators, the verification of the above properties reduces to checking  that the map $\si$ preserves
the ideals defining these algebras. In this sense we say that the defining relations of both algebras are compatible with the permutation map $\si$.

Note, that an obvious and well-known example of such a permutation map is the flip $\sigma = P$:
$$
P(a\otimes b) = b\otimes a\quad \forall\, a\in A,\,b\in B.
$$
If the permutation map $\sigma$ is not a flip (or a super flip), we call the corresponding double {\it a quantum double} (QD).

By means of the permutation map $\si$ we can introduce a unital associative algebra which is isomorphic as a vector space to the tensor product $B\ot A$. 
For this we consider a free tensor algebra $\mathrm{T}(A\oplus B)$ and take its quotient over an ideal
generated by the multiplication rules in $A$, $B$  and by elements $a\otimes b -\sigma(a\otimes b)$, $\forall \, a\in A$, $\forall\,b\in B$. So, in the quotient algebra the
following relations take place:
\be
a\otimes b = \sigma(a\otimes b).
\label{perm-rel}
\ee
Below the equalities (\ref{perm-rel}) will be referred to as {\it the permutation relations}.

Moreover, if the algebra $A$ is equipped with an algebra homomorphism $\varepsilon_A :A\to \C$, then it becomes possible to define a linear action 
$\triangleright: A\otimes B\rightarrow B$ of the algebra $A$ onto $B$:
\be
a\triangleright b:=(\mathrm{id}_B\ot \varepsilon_A) \si(a\ot b),\quad \forall\, \,a\in A,\, b\in B.
\label{act}
\ee
Note, that due to requirements imposed on the map $\sigma$ this action defines a representation of the algebra $A$ in the algebra $B$.

Consider now two examples of the QD which are of the main interest for us. As the first example we consider the QD $(A,B)$ of two RE algebras $A=\hLL(R)$ and
$B=\MM(R)$ equipped with the following permutation map
\be
\sigma:\,R_{12}\hL_1 R_{12}\otimes M_1\to M_1 \otimes R_{12}\hL_1 R^{-1}_{12}+R_{12}M_1\otimes 1_A.
\label{sii}
\ee
Consequently, the corresponding permutation relations are (from now on we omit the tensor product signs and the unit elements):
\be
R_{12}\hL_1 R_{12} M_1= M_1R_{12} \hL_1 R^{-1}_{12}+R_{12}M_1.
\label{dvaa}
\ee
This QD was introduced in \cite{GPS2, GPS3}. The modified RE algebra $A= \hLL(R)$ admits a homomorphism $\varepsilon: \hLL(R)\rightarrow {\Bbb C}$ which on generators
of $\hLL(R)$ is defined as $\varepsilon(\hL)=0$, $\varepsilon(1_A) = 1$. Then, in accordance with (\ref{act}) we get the action of the algebra $\hLL(R)$
onto $\MM(R)$:
$$
\hL_1 R_{12}\triangleright M_1=M_1.
$$
Note, that if we replace the Hecke symmetry $R$ to the flip $P$ in the permutation map (\ref{sii}), then the above action acquires the meaning of the left action of the enveloping
algebra\footnote{We assume the realization of $U(gl_N)$ as an algebra of right invariant differential opereators on $\mathrm{Sym}(gl_N)$.} $U(gl_N)$ onto the commutative algebra
$\mathrm{Sym}(gl_N)$.

Using the change (\ref{iso}) of the generating matrices, we can pass from the QD $(\hLL(R), \MM(R))$ to that $(\LL(R), \MM(R))$. The permutation relations in
the latter QD take the form
\be
R_{12}L_1R_{12} M_1 = M_1 R_{12} L_1R^{-1}_{12}.
\label{dva}
\ee
Since $\varepsilon(L)=I$ we get the action
\be
L_1R_{12}\triangleright M_1 = R^{-1}_{12}M_1.
\label{new-act}
\ee
From the technical point of view this form of the QD is easier to deal with, so we use it below for a study of invariant operators.

As the second example consider the QD $(\hLL(R), \MM(R))$ composed of the same algebras but equipped with another permutation map $\tilde\sigma$ leading
to the following permutation relations:
\be
R_{12}\hL_1R_{12} M_1-M_1 R_{12}\hL_1R_{12}= R_{12} M_1-M_1R_{12}.
\label{py}
\ee
With the homomorhism $\varepsilon(\hL)=0$ we get the action of the form:
$$
\hL_1R_{12}\triangleright M_1 = M_1-R^{-1}_{12}M_1R_{12}.
$$

Now, if we pass to the classical limit, assuming that  $R$ tends to $P$, the above action again turns into the action of the  enveloping algebra $U(gl_N)$ onto commutative
algebra $\mathrm{Sym}(gl_N)$ but the generators  of enveloping algebra should be realized as the adjoint vector fields.

Here, we can also express $\hL$ via $L$ by means of the relation (\ref{iso}). Then we get a QD $(\LL(R), \MM(R))$, equipped with the permutation relations of the form:
\be
R_{12}L_1R_{12} M_1= M_1R_{12} L_1 R_{12}.
\label{pyat}
\ee
As for the homomorphism $\varepsilon$, we always put $\varepsilon(L)=I$ for the generating matrix of the RE algebra $\LL(R)$ and $\varepsilon(\hL)=0$ for that of the modified 
RE algebra $\hLL(R)$.

Now, we go back to the QD $(\hLL(R), \MM(R))$ equipped with the permutation relations (\ref{dvaa}). As was shown in \cite{GPS3}, the  matrix $D=\|\pa_i^j\|$, defined by
$D= M^{-1}\hL$ satisfies the matrix relations:
$$
R^{-1}_{12}D_1 R^{-1}_{12}D_1=D_1 R^{-1}_{12}D_1 R^{-1}_{12}
$$
and
\be
D_1R_{12}M_1R_{12}= R_{12}M_1R^{-1}_{12}D_1+ R_{12}.
\label{raz}
\ee
The entries of the matrix $D$ generate an RE algebra\footnote{Observe that $R^{-1}$ is a braiding, subject to the Hecke condition with $q$ replaced by $\qq$.}  $\DD(R^{-1})$
which forms a QD with the algebra $\MM(R)$, equipped with the permutation relations (\ref{raz}). To define the action of $\DD(R^{-1})$ onto $\MM(R)$ we 
put $\varepsilon(D)=0$.

We treat the entries of the matrix $D$ as quantum analogs of the usual partial derivatives and call them quantum partial derivatives (QPD). This treatment is motivated by the fact
that if $R=P$ the generators $m_i^j$ become commutative and the entries of the matrix $D$ become the partial derivatives $\pa_i^j=\partial/\partial m_j^i$ while the permutation
relations (\ref{raz}) turn into the classical Leibniz rules. The QD $(\DD(R^{-1}), \MM(R))$ becomes the Heisenberg-Weyl algebra with the coordinates $(m_i^j)_{1\le i,j\le N}$ and the
corresponding partial derivatives $(\partial_i^j)_{1\le i,j\le N}$.

In the last section the QD $(\DD(R^{-1}), \MM(R))$ will be used for a definition of analogs of the partial derivatives on the enveloping algebra $U(gl_N)$.

We complete this section with one more example of a QD.  Let us take $A=\LL(R)$ again. As the algebra $B$ we choose the free tensor algebra $T(V)$ of the space $V$. Define
the permutation relations as follows:
\be
R_{12} L_1 R_{12} \,x_1=x_1L_2.
\label{perr}
\ee
On taking $\varepsilon(L) = I$, we get the action $L_1R_{12}\triangleright x_1=R_{12}^{-1}\,x_1$. Note that the permutation relations (\ref{perr}) can be reduced
to the $R$-symmetric and $R$-skew-symmetric algebras of $V$ defined respectively by
$$
\mathrm{Sym}_R(V)=T(V)/\langle \mathrm{Im}(q I-R)\rangle,\qquad \Lambda_R(V)=T(V)/\langle \mathrm{Im}(\qq I+R)\rangle.
$$

\section{Quantum left vector fields and the invariant operators}

Let $\hL$ be the generating matrix of the enveloping algebra $U(gl_N)$, that is the matrix $\hL$ meets the relation (\ref{mRE}) with $R=P$. It is well known that for all integers $k\ge 1$
the {\it power sums} $\mathrm{Tr} \hL^k$ belong to the center of $U(gl_N)$. Therefore, by the Schur lemma in any irreducible finite dimensional $U(gl_N)$-module $V_\lambda$ the
images of $\mathrm{Tr} \hL^k$ are scalar operators. Note that any such a module $V_\lambda$ is labelled by a partition $\la=(\la_1\geq \la_2\geq\dots \geq\la_N)$. The eigenvalues
of the mentioned scalar operators were computed in \cite{PP}.

Let us recall a way of constructing the $U(gl_N)$-modules $V_\lambda$. We start from the finite dimensional complex vector space $V$, $\dim_{\,\Bbb C}V = N$. In this space the
standard irreducible fundamental vector representation is realized.

Then for any integer $k\ge 1$ we consider the tensor power $V^{\otimes k}$. This space is decomposed into the direct sum of irreducible $U(gl_N)$-modules:
$$
V^{\otimes k} = \bigoplus_{\lambda\vdash k}\bigoplus_{T_\lambda}P_{T_\lambda}V^{\otimes k}.
$$
Here $T_\lambda$ denotes a standard Young table of the Young diagram, corresponding to a partition $\lambda\vdash k$. The orthonormal projection operators $P_{T_\lambda}$
are images of primitive idempotents $E_{T_\lambda}$ of the group algebra ${\Bbb C}[S_k]$ of the symmetric group $S_k$ under the representation of $S_k$ in $V^{\otimes k}$
where elements of the group $S_k$ act by permutations of factors of the tensor product $V^{\otimes k}$.

The subspaces $V_{T_\lambda} = P_{T_\lambda}V^{\otimes k}$ are irreducible $U(gl_N)$-modules. For different tables $T_\lambda$ and $T'_\lambda$ of the same diagram $\lambda$
these modules are isomorphic to each other and can be treated as different embeddings  of the module $V_\lambda$ into the space $V^{\otimes k}$:
$$
V_{T_\lambda}\simeq V_{T'_\lambda}\simeq V_\lambda.
$$

In what follows we assume the Hecke symmetry $R$ to be a deformation of the usual flip. In this case the corresponding modified RE algebra $\hLL(R)$ is a deformation of the
enveloping algebra $U(gl_N)$. Moreover, the category of finite dimensional $U(gl_N)$-modules can be also deformed into that of $\hLL(R)$-modules. Irreducible objects of this category can be constructed in the same way with the use of the Hecke algebras $H_k(q)$ instead of $\C[S_k]$. Namely, for a generic $q$ in the Hecke algebra $H_k(q)$ there are analogous idempotents $E_{T_\lambda}$ (we keep the same notation for them) such that $P_{T_\la}(R)V^{\ot k}$ are irreducible $\hLL(R)$-modules. Here
$P_{T_\lambda}(R) = \rho_R(E_{T_\lambda})$ and $\rho_R:H_k(q)\rightarrow \mathrm{End}(V^{\otimes k})$ is an $R$-matrix representation of the Hecke algebra in the space
$V^{\otimes k}$. This representation sends an Artin generator $\tau_i$ of the Hecke algebra to $R_i=I^{\otimes (i-1)}\otimes R\otimes I^{\otimes(k-i+1)}$. We refer the reader to
(\cite{GPS1}) for detail and to \cite{OP} for an explicit construction of the primitive idempotents $E_{T_\lambda}$ in the Hecke  algebra.  

An important role in the theory of the Hecke algebras belongs to the so-called Jusys-Murphy elements which generates the maximal commutative subagebra in the Hecke algebra
$H_k(q)$. In the $R$-matrix representation they form a set of mutually commutative linear operators $J_i$, $1\le i\le k$, in the space $V^{\otimes k}$:
$$
J_1=\mathrm{Id}_{V^{\otimes k}},\quad  J_i=R_{i-1}J_{i-1}R_{i-1},\quad 2\leq i\leq k.
$$

Introduce also the notation for some ``matrix copies'' of the generating matrix $L$:
$$
L_{\ov 1}= L_{\un 1}=L_{1},\quad L_{\ov k}=R_{k-1}L_{\ov {k-1}}\,R^{-1}_{k-1},\quad L_{\un k}=R^{-1}_{k-1}L_{\un {k-1}}R_{k-1}.
$$

With the above notation we can write the action (\ref{new-act}) in the form:
\be
L_{\un 2}\triangleright M_1=J_2^{-1} M_1.
\label{JM}
\ee
On applying the permutation relations (\ref{dva}) one can generalize this formula on arbitraty ``monomials'' in generators of $\MM(R)$:
\be
L_{\un{k+1}}\triangleright M_1\dots M_{\ov k}=J_{k+1}^{-1} M_1\dots M_{\ov k}.
\label{dvaaa}
\ee

Calculating the $R$-trace at the $(k+1)$-th space and taking into account the formula
\be
\Tr_{R(k+1)} X_{\ov {k+1}}=\Tr_{R(k+1)} X_{\un {k+1}}= I_{1\dots k} \Tr_{R} \,X,
\label{syst}
\ee
where $X$ is an arbitrary $N\times N$ matrix, we get the following claim (see \cite{S}).

\begin{proposition} \label{prop:2}
Let $\lambda\vdash k$ be a partition and $T_\lambda$ be a standard Young table corresponding to the Young diagram of the partition $\lambda$.
Then the following relation holds
\be
\Tr_R L\triangleright P_{T_\lambda}(R) M_1M_{\ov 2}\dots M_{\ov k} =\chi_{\la} (\Tr_R L)P_{T_\lambda}(R)M_1M_{\ov 2}\dots M_{\ov k},  \label{spec}
\ee
where
\be
\chi_{\la} (\Trr L)=\frac{N_q}{q^{N}}-\frac{\nu}{q^{2N}} \sum_{i=1}^k  \, q^{-2c(i)},\quad \nu = q-q^{-1},
\label{dec}
\ee
and the sum is taken over all boxes of the table $T_\lambda$. Here  $c(i)=n-m$ is the content of the box in which the integer $i$ is located, that is $n$ (respectively $m$) is the 
number of the corresponding column (row).
\end{proposition}

Note that the symbol of the table $T$ in the notation $\chi_\la$  would be irrelevant since the right hand side of (\ref{dec}) actually depens only on the {\it  diagram} $\lambda$ but not
on the table.

\begin{remark} 
Let us point out  that formula (\ref{dvaaa}) is the matrix form of the action of the left quantum  vector fields on the algebra $\MM(R)$. In fact, it is similar to the action of the algebra
$\LL(R)$ onto the free tensor algebra $T(V)$ considered at the end of the previous section. Namely, using the permutation relations (\ref{perr}) and the homomorphism
$\varepsilon(L) = I$ we can get the action
\be
L_{\un{k+1}}\triangleright x_1\ot\dots \ot x_{ k}=J_{k+1}^{-1}x_1\ot\dots \ot x_{ k}.
\label{tri}
\ee
Evidently, a formula similar to (\ref{spec}) is also valid in this case.
\end{remark}

We are interested in the spectral analysis of the invariant operators, that is these corresponding to the central elements of the algebra $\LL(R)$ acting onto $\MM(R)$. The proposition
\ref{prop:2} describes the spectrum of the operator $\Trr L$, corresponding to the lowest power sum.

Below,  we consider some other examles of central elements and describe the spectrum of the corresponding operators in invariant subspaces $\mathrm{Im}(P_{T_\lambda})$ extracted
by the projection operators $P_{T_\lambda}(R)$. Besides the higher power sums $\mathrm{Tr}_RL^k$, $k\ge 2$, we are interested in the so-called quantum elementary symmetric polynomials $e_k(L)$, which are defined by the formula:
$$
e_k(L)=\Tr_{R(1\dots k)} (A^{(k)} L_1 L_{\ov 2}\dots L_{\ov k}),\quad  k \ge 1.
$$
Here $A^{(k)}$ are the $R$-skew-symmetrizers introduced by the following recursion formula:
$$
A^{(1)}=I,\quad  A^{(k)}=\frac{1}{k_q}A^{(k-1)}\left(q^{k-1} I-(k-1)_q \, R_{k-1}\right)A^{(k-1)},\quad k\ge 2,
$$
where $k_q=(q^k-q^{-k})/(q-\qq)$ is a $q$-integer. Below we assume $q^{2k}\not=1$ for all integer $2\le k\le N$, and, therefore, $k_q\not=0$. Any $R$-skew-symmetrizer $A^{(k)}(R)$ is a particular case of the above idempotents $P_{T_\lambda}(R)$, corresponding to one-column diagram $\la=(1^k)$.

The elementary symmetric polynomials possess a natural and convenient parameterization in terms of the {\it quantum eigenvalues} $\{\mu_i\}_{1\le i\le N}$ of the generating
matrix $L$. These eigenvalues belong to an algebraic extension of the center of the RE algebra $\LL(R)$. They are defined via the quantum Cay\-ley-Ha\-mil\-ton
identity on the matrix $L$:
\be
L^N-q e_1(L) L^{N-1}+ q^2e_2(L) L^{N-2}+\dots  +(-q)^{N} e_{N}(L) I=0.
\label{CH}
\ee
Let us introduce $N$ elements $\mu_i$, $1\le i\le N$, which satisfy a system of $N$ polynomial equations with central elements on the left hand side:
\be
 q^k e_k(L)=\sum_{1\leq i_1<\dots <i_k\leq N} \mu_{i_1}\dots \mu_{i_k}, \quad 1\le k\le N.
\label{elem-mu}
\ee
The  elements $\mu_i$ are also assumed to be central. They are interpreted as quantum analogs of eigenvalues of the matrix $L$ since in virtue of (\ref{elem-mu}) the 
Cayley-Hamilton identity (\ref{CH}) can be rewritten in the following factorized form:
\be
\prod_{i=1}^N(L-\mu_iI) = 0.
\label{factor-CH}
\ee
In terms of the quantum eigenvalues $\mu_i$ one can express all central symmetric polynomials of the RE algebra --- the power sums, the quantum Schur functions and so on
(see \cite{GPS1,GPS4} for more detail). For example, the power sums in terms of quantum spectrum read:
\be
\Trr L^k=\sum \mu_i^k \, d_i,\qquad d_i=\qq \prod_{j\not=i}^N \frac{\mu_i-q^{-2}\mu_j}{\mu_i-\mu_j}.
\label{powe}
\ee
Since the elements $\mu_i$ are central, then in any irreducible $\LL(R)$-module they are represented by scalar operators. These operators must be compatible with relations 
(\ref{elem-mu}).  In the paper \cite{GPeS1} we suggested the following conjecture\footnote{In the cited paper the conjecture was formulated for a wider class of Hecke symmetries, 
namely those of the bi-rank $(m|0)$. (For the notion of the bi-rank the reader is referred to \cite{GPS1}.) The Hecke symmetries we are dealing with in the present paper is of the 
bi-rank $(N|0)$.}.
\begin{conjecture}
In $\LL(R)$-invariant subspaces $V_{T_\lambda} = \mathrm{Im}P_{T_\lambda}(R)$ the central elements $\mu_i$ $1\le i\le N$ are represented (after a proper ordering) by the 
following scalar operators:
\be
\mu_i\mapsto \chi_\lambda(\mu_i)\,\mathrm{Id}_{V_{T_\lambda}},\qquad \chi_\la(\mu_i)=q^{-2(\la_i+N-i)},
\label{mu-char}
\ee
where $\lambda_i$ are elements of a partition $\lambda = (\lambda_1,\lambda_2,\dots ,\lambda_N)$.

The representation (\ref{mu-char}) is compatible with the scalar operators corresponding to the quantum elementary symmetric polynomials 
$e_k(L)\mapsto \chi_\lambda(e_k)\,\mathrm{Id}_{V_{T_\lambda}}$, that is:
$$
q^k\,\chi_{\la}(e_k)=\sum_{1\leq i_{1}<\dots < i_k \leq N}   \chi_{\la}(\mu_{i_1})\dots \chi_{\la}(\mu_{i_k}).
$$
\end{conjecture}
For the lowest elementary symmetric polynomial $e_1(L) = \mathrm{Tr}_R(L) = q^{-1}\sum_i\mu_i$ this conjecture is easily verified with the use of (\ref{dec}).
We have also checked it for $e_2(L)$ by a direct computation.

Turn now to the modified RE algebra $\hLL(R)$ whose generators $\hat l_i^j$ are treated as quantum analogs of the right-invariant vector fields. The corresponding Casimir 
operators are $\Trr \hL^k$. The most convenient way to perform their spectral analysis is to introduce the quantum eigenvalues $\hat\mu_i$ of the generating matrix $\hat L$. 
Taking into account relation (\ref{iso}), we can rewrite the factorized Cayley-Hamilton identity (\ref{factor-CH}) in the form of the matrix identity for $\hat L$:
$$
\prod_{i=1}^N(L-\mu_i I) = 0\qquad \Rightarrow\qquad \prod_{i=1}^N(\hat L-\hat \mu_i I) = 0,
$$
where the quantum eigenvalues $\hat \mu_i$ of the generating matrix $\hat L$ are connected with those of $L$ by a linear shift:
\be
\mu_i=1-\nu \hat \mu_i,\qquad \nu = q-q^{-1}.
\label{mu-shift}
\ee
Here the symbol $1$ stands for the unit element of the RE algebra.

As a direct consequence of (\ref{mu-shift}) we obtain:
$$
\chi_\la({\hat{\mu}_k})=\frac{1- q^{-2\, (\la_k+N-k)}}{q-\qq} = q^{-(\lambda_k+N-k)}(\lambda_k+N-k)_q,\quad 1\le k \le N.
$$

By taking the limit $q\rightarrow1$, we get the spectral values of $U(gl_N)$ generating matrix $\hL$ in the irreducible module $V_\lambda$:
$$
\chi_\la({\hat{\mu}_k})=\la_k+N-k,\quad 1\le k\le N.
$$
Note that this is in the full agreement with the results of \cite{PP}.

\section{Different forms of the Capelli identity}

Consider the QD $(\DD(R^{-1}), \MM(R))$ of two RE algebras generated by the enties of matrices $D$ and $M$ and equipped with the permutation relations (\ref{raz}).
Introduce the matrix $\hat L=MD$. As was argued in Section 2, its entries generate the modified RE algebra $\hLL(R)$ and the permutation relations with entries of the matrix $M$
are given by (\ref{dvaa}).

One of the remarkable properties of the QD under consideration are the following quantum {\it matrix Capelli identities}\footnote{If $R$ is a deformation of the flip $P$ then at the
classical limit $q\rightarrow 1$ these identities transform into those presented in \cite{Ok2}.}.
\begin{theorem}  {\rm (\cite{GPeS2}) }
\label{th:1} 
For $\forall\,k\ge 1$ the following matrix identity takes place:
\be
A^{(k)}\hat L_{\overline 1}\,(\hat L_{\ov 2}\,+q I)\dots (\hat L_{\ov k}\,+q^{k-1}(k-1)_qI\,) \,A^{(k)}=
q^{k(k-1)}A^{(k)}M_{\overline 1}\dots M_{\ov k}\, D_{\ov k}\dots D_{\overline 1}.
\label{th}
\ee
\end{theorem}

Let us point out that identities (\ref{th}) are valid for a wide class of the Hecke symmetries: to prove (\ref{th}) one needs only the braid relation and the Hecke condition (\ref{H-cond})
on $R$. If we impose some additional restrictions on $R$ we can obtain further consequences of (\ref{th}), namely, the quantum version of the determinant Capelli identity.

Thus, assume $R$ to be a deformation of the usual flip $P$, then the $R$-skew-symmetrizer $A^{(N)}$ is a unit rank projector
\be
\mathrm{dim}\,\mathrm{Im}\, A^{(N)}=1
\label{unit-rank}
\ee
and $A^{(N+1)}\equiv 0$ (recall that $N=\dim\, V$). A direct consequence of (\ref{unit-rank}) is the existence of two tensors $|u\rangle= \|u_{i_1i_2\dots \,i_N}\|$ 
and $\langle v| = \|v^{i_1i_2\dots\,i_N}\|$ such that
$$
A^{(N)}  = |u\rangle\otimes\langle v|,  \qquad \langle v| u\rangle = \sum_{\{i\}} v^{i_1i_2\dots\,i_N}u_{i_1i_2\dots\,i_N} = 1.
$$
The structure tensors $|u\rangle$ and $\langle v|$ allow one to define a quntum analog of the determinant for the generatig matrix of the RE algebra. Namely, we set
$$
{\det}_R\,M =\langle v| M_1M_{\overline 2}\dots M_{\overline N} |u\rangle,\qquad {\det}_{R^{-1}}\,D = \langle v| D_{\ov N}\dots D_{\overline 2}D_1)|u \rangle .
$$
With these definitions we can formulate the following corollary of Theorem \ref{th:1}.
\begin{corollary} \label{cor:6}
If $R$ is a deformation of the flip $P$, then the following quantum Capelli identity takes place:
\be
\mathrm{Tr}_{R(1\dots N)}A^{(N)}\hat L_{\overline 1}\,(\hat L_{\ov 2}\,+q I)\dots (\hat L_{\ov N}\,+q^{N-1}(N-1)_qI\,) =
q^{-N}{\det}_R\, M\, {\det}_{R^{-1}}\, D.
\label{dettt}
\ee
\end{corollary}
\begin{remark}
Note that this corollary remains valid mutatis mutandis if the initial Hecke symmetry $R$ is of bi-rank $(m|0)$, $m\leq N$, the corresponding examples were constructed in \cite{G}. 
This condition means that $\mathrm{rank}\,A^{(m)}=1$ and $A^{(m+1)} \equiv 0$. In this case we only have to replace $N$ by $m$ in formula (\ref{dettt}), while the structure tensors 
$|u\rangle$ and $\langle v|$ should be extracted from the projector $A^{(m)}$.
\end{remark}

In conclusion of the section we note that another version of the quantum Capelli identity was suggested in \cite{NUW}. This version is related to the RTT algebra and is valid for $R$
coming from the QG $U_q(sl_N)$.

\section{Quantum adjoint vector fields and quantum orbits}

In this section we deal with the QD $(\hLL(R), \MM(R))$, endowed with the permutation relations (\ref{py}). Our main objective here is to show that the action of the modified RE algebra
$\hLL(R)$ on $\MM(R)$ can be reduced onto some quotients of the algebra $\MM(R)$. These quotient algebras are treated as quantum analogs of the coordinate rings of the 
$GL(N)$-orbits in $gl_N^*$. Observe that we employ the term ``orbit"  in a loose sense since even in the classical case these quotients are ``orbits" for generic matrices\footnote{This 
means that the eigenvalues of the matrix are pairwise distinct.} generating the orbit. The question which quantum orbits can be considered as ``generic" is discussed below. The 
generators $\hat l_i^j$ are treated as quantum adjoint vector fields.

\begin{proposition}\label{prop:8}
In the QD $(\hLL(R),  \MM(R))$ the following relations hold for any integer $k\ge 1$ 
$$
\hL\, \Trr( M^k)=\Trr( M^k)\hL
$$
and consequently
\be
\hL\triangleright \Trr M^k=0.
\label{hL-act}
\ee
\end{proposition}

\medskip

\noindent
{\bf Proof.}
Successively multiplying the both sides of matrix equality (\ref{py}) by the matrix $M_1$ from the right and transforming the term $M_1R_{12}\hat L_1 R_{12}M_1$
with the use of (\ref{py}) we come to the result:
$$
R_{12}\hat L_1R_{12} M^k_1 - M_1^kR_{12}\hat L_1 R_{12} = R_{12}M_1^k - M_1^kR_{12}.
$$
Then we multiply this equality by $R^{-1}_{12}$ from the left and right and apply the $R$-trace in the second space. These operations lead to the desired result
in virtue of (\ref{syst}):
$$
\hL_1 \Trr M^k - \Trr M^k\hL_1 = 0.
$$
To get the action (\ref{hL-act}) we apply the homomorphism $\varepsilon(\hat L) = 0$ to the second term.\hfill\rule{6.5pt}{6.5pt}

\medskip

So, the quantum adjoint vector fields kill all elements $\Trr M^k$. Note that these elements are central in the RE algebra $\MM(R)$.

Now, consider the following quotient algebra:
$$
\mathcal{O}(\al_1,\dots ,\al_N)=\MM(R)/\langle \Trr M-\al_1, \, \Trr M^2-\al_2,\dots , \Trr M^N-\al_N\rangle,
$$
where $\al_1,\dots ,\al_N$ are some complex constants and the notation $\langle J\rangle$ stands for the ideal, generated by a subset $J$ of a given algebra.
We call these quotients the ``quantum orbits".

The quantities $\al_i$ can be expressed in terms of the eigenvalues $\mu_i$ of the matrix $M$ since for this matrix formula (\ref{powe}) is also valid. In the paper \cite{GS1} we
showed (up to a conjecture on the quantum de Rham complex) that the quantum orbits are generic if
\be
\mu_i\not=q^2\mu_j,\quad\forall \, i,j.
\label{orb}
\ee
Observe that by definition a quantum orbit is generic if its cotangent module is projective. Thus, as was shown in \cite{GS1}, under the condition (\ref{orb}) and under the assumption
that the mentioned conjecture is true the cotangent module on the quantum orbit $\mathcal{O}(\al_1,\dots ,\al_N)$ is projective.

\section{Quantum partial derivatives on $U(gl_N)$ background}

In this section we consider one more application of the quantum left vector fields. Namely, we show how to get analogs of the partial derivatives on the enveloping algebra $U(gl_N)$
via the QD introduced at the end of Section 2.

To this purpose we return to the QD $(\DD(R^{-1}), \MM(R))$, endowed with the permutation relations (\ref{raz}). Let us introduce a new generating matrix $N$ of the RE algebra $\MM(R)$ by the shift analogous to (\ref{iso}):
$$
M=h\, I-(q-\qq)\, N,\qquad h\in{\Bbb C}.
$$
The full list of the relations in the QD $(\DD(R^{-1}), \hat\NN_h(R))$, generated by the entries of matrices $D$ and $N$, is as follows
$$
R_{12} N_1R_{12} N_1-N_1R_{12} N_1R_{12}=h(R_{12} N_1-N_1 R_{12}),
$$
$$
 R^{-1}_{12} D_1 R^{-1}_{12} D_1=D_1R^{-1}_{12} D_1 R^{-1}_{12},
$$
\be
D_1R_{12}N_1R_{12}- R_{12}N_1R^{-1}_{12}D_1=R_{12}+h D_1R_{12}.
\label{q-Leib-rule}
\ee
We introduce the deformation parameter $\h$ in order to present the modified RE algebra $\hat\NN_h(R)$ as a deformation of the RE algebra $\MM(R)$, which corresponds to the
value $\h=0$.

The entries of the matrix $D=\|\pa_i^j\|_{1\leq i,j \leq N}$ can be treated as analogs of partial derivatives in $n_i^j$ and we call them the quantum partial derivatives (QPD). Using the
homomorphism $\varepsilon(D) = 0$, we get the  action of the QPD on the generators of $\hat\NN_h(R)$:
\be
D_1\triangleright N_{\overline 2} = R_{12}^{-1}.
\label{lin-act}
\ee
The relation (\ref{q-Leib-rule}) is a substitution of the Leibniz rule allowing one to extend the action of QPD (\ref{lin-act}) from generators onto the whole algebra $\hat\NN_h(R)$.

If $R$ is a deformation of the flip $P$ then at the limit $q\rightarrow 1$ relations (\ref{lin-act}) turn into the classical action of the partial derivatives on the generators of the commutative
algebra\footnote{However, the action of the QPD on higher polynomials in $n_k^l$ differs from the classical formulae and contains corrections depending on $h$.
The full classical picture restores after specializing $h=0$.} $\mathrm{Sym}(gl_N)$:
$$
\partial_i^j\triangleright n_k^s = \delta_i^s\,\delta_k^j.
$$
This is the reason for treating the QPD $\partial_i^j$ as ``derivatives'' in generators $n_j^i$.

In general, if the symmetry $R$ is involutive, it is useful to introduce the matrix $\hD=D+\h^{-1}\, I$. In terms of this matrix the permutation relations (\ref{q-Leib-rule}) take the form:
$$
\hD_1R_{12}N_1R_{12} - R_{12}N_1R_{12}\hD_1= h \hD_1R_{12}
$$
We call the entries of the matrix $\hD=\|\hpa_i^j\|$ the shifted QPD\footnote{ In fact, only the diagonal elements of the matrix $D$ are shifted.}.

From now on we assume that $R=P$. In this case the permutation relations for the generators of $\NN_h(R)$ turn into those of the universal enveloping algebra $U(gl(N)_h)$:
$$
N_2N_1-N_1N_2 = h(P_{12}N_1 - N_1P_{12}).
$$
The permutation relations (\ref{q-Leib-rule}) give rise to the Leibnitz rule for the QPD, expressed via the coproduct defined on the shifted QPD
$$
\De(\hpa_i^j)=\hpa_k^j\ot \hpa_i^k,
$$
and completed with $\varepsilon(\hpa_i^j)=h^{-1}\, \de_i^j$.

Let us consider the example $N=2$ in more detail. It is convenient to pass to the compact form $u(2)_h$ of the algebra $gl(2)_h$. The Lie algebra $u(2)_h$ is generated  by 4 elements
$x$, $y$, $z$ and $t$ with the following Lie brackets:
$$
[x, \, y]=h z,\qquad [y, \, z]=h x,\qquad[z, \, x]=h y,\qquad [t, \, x]=[t, \, y]=[t, \, z]=0.
$$

The corresponding QPD $\pa_x$, $\pa_y$, $\pa_z$ and $\hat\pa_t$ (we use the ``shifted" derivative in  $t$:  $\hpa_t = \pa_t + 2/h$) commute with each other,
while the permutation relations with the $u(2)_h$ generators read:
\be
\begin{array}{l@{\quad}l@{\quad}l@{\quad}l}
\hpa_t\,t - t\,\hpa_t = \frac{h}{2}\,\hpa_t & \hpa_t\, x - x\,\hpa_t
=-\frac{h}{2}\,\pa_x &
\hpa_t\, y - y\, \hpa_t=-\frac{h}{2}\,\pa_y &\hpa_t\, z - z\,\hpa_t=- \frac{h}{2}\,\pa_z\\
\rule{0pt}{7mm}
\pa_x\, t - t\,\pa_x = \frac{h}{2}\,\pa_x &\pa_x \,x -  x\,\pa_x = \frac{h}{2}\,\hpa_t &
\pa_x \, y-  y\,\pa_x = \frac{h}{2}\,\pa_z & \pa_x \,z - z\, \pa_x  = - \frac{h}{2}\,\pa_y \\
\rule{0pt}{7mm}
\pa_y \,t - t \, \pa_y = \frac{h}{2}\,\pa_y & \pa_y \,x -  x\,  \pa_y = -\frac{h}{2}\,\pa_z &
\pa_y \,y - y \,  \pa_y = \frac{h}{2}\,\hpa_t & \pa_y \,z - z \,  \pa_ y= \frac{h}{2}\,\pa_x\\
\rule{0pt}{7mm}
\pa_z \,t - t \,\pa_z = \frac{h}{2}\,\pa_z & \pa_z \,x - x \,\pa_z = \frac{h}{2}\,\pa_y&
\pa_z \,y -  y\,\pa_z = -\frac{h}{2}\,\pa_x & \pa_z \,z - z \,\pa_z = \frac{h}{2}\,\hpa_t.
\end{array}
\label{leib}
\ee

Introduce the matrix $\hatt$, composed of the partial derivatives:
\be
{\hatt}=i\hh \left(\begin{array}{rrrr}
\hpa_t&\pa_x&\pa_y&\pa_z\\
-\pa_x&\hpa_t&-\pa_z&\pa_y\\
-\pa_y&\pa_z&\hpa_t&-\pa_x\\
-\pa_z&-\pa_y&\pa_x&\hpa_t
\end{array} \right),
\label{seven}
\ee
where $\hh=h/2i$.

By $\hatt(a)$ we denote the matrix whose entries result from applying the corresponding partial derivatives to an element $a\in U(u(2)_h)$. Thus, we get a linear map of associative
algebras:
$$
U(u(2)_h)\to \mathrm{Mat}_4(U(u(2)_h)) = \mathrm{End}({\Bbb C}^{\,4})\otimes U(u(2)_h)
$$
which will be also denoted by the  symbol $\hatt$:
\be
\hatt: \, a\mapsto  \hatt(a)\quad \, \forall\, a\in U(u(2)_h).
\label{mapp}
\ee

We state (see \cite{GS3}) that the map $\hatt$  is actually a homomorphism of associative algebras, that is
\be
\hatt(ab)=\hatt(a) \cdot \hatt(b),\quad \forall\,a,b\in U(u(2)_h).
\label{mult}
\ee

Consequently, the map $\hatt$ defines a  representation of the algebra $U(u(2)_h)$:
$$
\hatt(x)\cdot\hatt(y)-\hatt(y)\cdot\hatt(x)=h\hatt(z)
$$
and so on. In virtue of the above properties the map $\hatt$ is treated as the Leibniz rule for the QPD on the algebra $U(u(2)_h)$.

We want to stress the usefulness of this form of the Leibniz rule. So far, the QPD were defined on the algebra of polynomials in noncommutative generators. Our next aim is to
extend their action onto some elements of a central extension of this algebra. Namely, consider the so-called {\em quantum radius} $\rh$:
$$
\rh^2=x^2+y^2+z^2+\hh^2,
$$
which is an element of a central extension of the algebra $U(u(2)_h)$. We want to extend the action of the QPD onto $\rh$. As a reasonable criterium of such an extension we require
that the property (\ref{mult}) would be preserved on the extended algebra. Namely, we require the relation
$$
\hatt(\rh)^2=\hatt(\rh^2)=\hatt(x^2+y^2+z^2+\hh^2)
$$
to be fulfilled.

As follows from computation of \cite{GS3}, the result of applying the map $\hatt$ to the quantum radius $\rh$ is
$$
{\hatt}(\rh)=\frac{\rh^2+\hh^2}{\rh}\, I+\frac{i\hh}{\rh}\, M,
$$
where the matrix $M$ is of the form:
$$
M=
\left(\!\!
\begin{array}{cccc}
 0 &    x&y&z\\
 -x&0&-z&y\\
-y&z&  0 & -x\\
-z&-y& x & 0
\end{array}\right).
$$
Now we are able to compute the action of the QPD on the quantum radius:
$$
\pa_t \rh=-\frac{i\hh}{\rh},\qquad \pa_x\rh=\frac{x}{\rh},\qquad \pa_y\rh=\frac{y}{\rh},\qquad \pa_z\rh=\frac{z}{\rh}.
$$
Observe that these formulae possess the usual classical limit as $\hh\to 0$ and, consequently, $\rh\to r$.

Also, in \cite{GS4} the following problem was discussed: how to extend the action of the QPD on the skew-field $B[B^{-1}]$, where $B$ is a central extension of the algebra
$U(u(2)_h)$. Assuming that the above Leibniz rule is still valid on this skew-field, we have:
$$
\hatt (a)\hatt (a^{-1})=\hatt (e)=I.
$$
Thus, in order to extend the action of the partial derivatives on the skew-field $B[B^{-1}]$ it suffices to invert the matrix $\hatt (a)$. In \cite{GS4} we exhibited an example,
in which such inverting was performed with the use of the Cayley-Hamilton identity for the matrix $\hatt (a)$.

In conclusion, we note that defining QPD on the extended algebra $U(u(2)_h)$ enables us to introduce some dynamical models on this algebra. Certain of them were considered in
\cite{GS2, GS3}.


\begin{thebibliography}{GPeS0}
\bibitem[G]{G} D. Gurevich,  Algebraic aspects of quantum Yang-Baxter equation, Leningrad Math. J., 2 (1990) 119--148.

\bibitem[GPeS1]{GPeS1} D. Gurevich, V. Petrova, P. Saponov,  $q$-Casimir and $q$-cut-and-join operators related to Reflection equation algebras,  Arxiv:2110.04354.

\bibitem[GPeS2]{GPeS2} D. Gurevich, V. Petrova, P. Saponov, Matrix Capelli identities related to reflection equation algebra, J.Geometry and Physics, 179 (2022) 104606.

\bibitem[GPS1]{GPS1} D. Gurevich, P. Pyatov, P. Saponov, Representation theory of (modified) Reflection Equation Algebra of the $GL(m|n)$ type,
St Petersburg Math. J., 20 (2009) 213--253.

\bibitem[GPS2]{GPS2} D. Gurevich, P. Pyatov, P. Saponov, Braided Differential operators on quantum algebras, J. of Geometry and Physics, 61 (2011) 1485--1501.

\bibitem[GPS3]{GPS3} D. Gurevich, P. Pyatov, P. Saponov, Braided Weyl algebras and differential calculus on $U(u(2))$, J. of Geometry and Physics, 62 (2012) 1175--1188.

\bibitem[GPS4]{GPS4} D. Gurevich, P. Pyatov, P. Saponov, Spectral parametrization for power sums of quantum supermatrices, Theor. and Math. Physics, 159 (2009) 587--597.

\bibitem[GS1]{GS1} D. Gurevich, P. Saponov, Geometry of non-commutative orbits related to Hecke symmetries, Contemporary Mathematics, 433 (2007) 209--250.

\bibitem[GS2]{GS2} D. Gurevich, P. Saponov, Noncommutative Geometry and dynamical models on $U(u(2))$ background, Journal of Generalized Lie theory and Applications, 
9:1 (2015).

\bibitem[GS3]{GS3} D. Gurevich, P. Saponov, Quantum geometry and quantization on $U(u(2))$ background. Noncommutative Dirac monopole, J. of Geometry and Physics, 
106 (2016) 87--97.

\bibitem[GS4]{GS4} D. Gurevich, P. Saponov, Noncommutative geometry on central extension of $U(gl(2))$, https://arxiv.org/abs/2009.05807.

\bibitem[IP]{IP} A. Isaev A., P. Pyatov, Spectral extension of the quantum group cotangent bundle, Comm. Math. Phys., 288 (2009) 1137--1179.

\bibitem[NUW]{NUW} M. Noumi, T. Umeda, M. Wakayama, A quantum analogue of the Capelli identity and elementary differential calculus  on $GL_q(n)$, Duke Math.J.,
76 (1994) 567--594.

\bibitem[OP]{OP} O. Ogievetsky, P. Pyatov, Lecture on Hecke algebras, In “Symmetries and Integrable systems”, Dubna publishing 2000,
Preprint CPT-2000/P.4076.

\bibitem[Ok1]{Ok1} A. Okounkov, Quantum immanants and higher Capelli identities, Transformation groups, 1 (1996) 99--126.

\bibitem[Ok2]{Ok2} A. Okounkov, Young basis, Wick formula and higher Capelli identities, International Mathematics Research Notices, 17 (1996) 817--839.

\bibitem[PP]{PP} A. Perelomov, V. Popov, Casimir operators for the classical groups, Dokl. Akad. Nauk SSSR, 174 (1967) 287--290.

\bibitem[S]{S} P. Saponov, Weyl approach to representation theory of reflection equation algebra, Journal of Physics A: Mathematical and General,
 vol. 37, no. 18 (2004) 5021--5046.

\end{thebibliography}
\end{document}